# RELATIVE FREQUENCIES IN MULTITYPE BRANCHING PROCESSES*

By Andrei Y. Yakovlev[1] and Nikolay M. Yanev[1,2]

*University of Rochester and Bulgarian Academy of Sciences*

This paper considers the relative frequencies of distinct types of individuals in multitype branching processes. We prove that the frequencies are asymptotically multivariate normal when the initial number of ancestors is large and the time of observation is fixed. The result is valid for any branching process with a finite number of types; the only assumption required is that of independent individual evolutions. The problem under consideration is motivated by applications in the area of cell biology. Specifically, the reported limiting results are of advantage in cell kinetics studies where the relative frequencies but not the absolute cell counts are accessible to measurement. Relevant statistical applications are discussed in the context of asymptotic maximum likelihood inference for multitype branching processes.

**1. Introduction.** The paper deals with multitype branching processes with discrete or continuous time assuming only the usual independence of the individual evolutions. Such processes are considered from a new perspective: modeling of the relative frequencies (proportions, fractions) of different types of individuals (instead of their absolute numbers) as functions of time. Notwithstanding other possible applications, we will use the term "cell" instead of referring to an abstract individual or particle throughout this paper to emphasize its special relevance to cell biology.

It is well known (see, e.g., Athreya and Ney [1] and Mode [19]) that in the positively regular and supercritical case the multitype population vector

Received May 2007; revised February 2008.
[1]Supported by NIH/NINDS Grant NS-39511 and NIH Grant N01-AI-050020.
[2]Supported in part by ECO-NET-06 and Grant VU-MI-105/2005.
[*]Shortly after the final version was submitted, Professor Andrei Yakovlev died on February 27, 2008. This paper is a tribute to the stimulating ideas and friendship that he shared in this collaboration.
*AMS 2000 subject classifications.* Primary 60J80, 60J85; secondary 62P10, 92D25.
*Key words and phrases.* Multitype branching processes, relative frequencies, joint asymptotic normality, cell proliferation.

$\mathbf{Z}(t)$ is asymptotically proportional (on the nonextinction set) to the left-eigenvector $v = (v_1, v_2, \ldots, v_d)$ corresponding to the Perron–Frobenius root of the mean matrix. Then the frequency of the $k$th type converges a.s. to $v_k / \sum_{j=1}^{d} v_j$ (see Mode [19], Theorem 8.3, Chapter 1). Jagers [15] was probably the first to consider relative frequencies in the context of biological applications. He studied asymptotic properties (as $t \to \infty$) of a reducible age-dependent branching process with two types of cells and proved convergence of their relative frequencies to nonrandom limits in mean square and almost surely on the nonextinction set. The usefulness of such frequencies in cell cycle analysis was further demonstrated by Mode [19, 20] who considered a four-type irreducible age-dependent branching process. Some of these results are discussed in Section 6. The relative frequencies of distinct cell types play an important role in the analysis of biological studies of proliferation and differentiation of cells. As discussed in Section 2, the need for such characteristics of cell kinetics arises in experimental situations where absolute cell counts are not accessible to measurement.

To formulate the problem in more technical terms, consider a multitype branching stochastic process $\mathbf{Z}(t) = (Z_1(t), Z_2(t), \ldots, Z_d(t))$, where $Z_k(t)$ denotes the number of cells of type $T_k$ $(k = 1, 2, \ldots, d)$ at time $t$. The time may be discrete $(t \in \mathbf{N}_0 = \{0, 1, 2, \ldots\})$ or continuous $(t \in \mathbf{R}^+ = [0, \infty))$, and the evolution of cells may be arbitrarily complex as long as it can be described by a (decomposable or nondecomposable) branching process. The main focus in the theory of branching process is on probabilistic characteristics of the process $\mathbf{Z}(t)$ and their asymptotic behavior when $t \to \infty$. In practical applications, however, the investigator is more interested in such characteristics at fixed time points (transient process) that provide typically much more information on model parameters than their asymptotic counterparts. It is still reasonable and useful to study another asymptotic aspect of the problem, namely, the behavior of $\mathbf{Z}(t)$ when the initial number of cells is large, which situation is of frequent occurrence in cell kinetics studies. This aspect becomes especially important when considering relative frequencies of cell types rather than the process $\mathbf{Z}(t)$ itself.

We define relative frequencies of cell types as $\Delta_k(t) = Z_k(t)/U(t), k = 1, 2, \ldots, d$, for $U(t) > 0$, where $U(t) = \sum_{k=1}^{d} Z_k(t)$ is the total number of cells present in the population at time $t$. In what follows, we assume that there is at least one type (say, type $T_1$) giving rise to all other types, regardless of whether the process is decomposable or nondecomposable. We assume in addition that the process begins with $Z_1(0) = N$ cells and study asymptotic properties of the above-defined fractions $\Delta_k(t; N)$ as $N \to \infty$.

The paper is organized as follows. Section 2 provides a biological motivation for the mathematical problem under consideration. The basic notions and preliminaries are introduced in Section 3. Section 4 includes Theorem 1 stating multivariate asymptotic normality of the fractions $(\Delta_1(t; N), \Delta_2(t; N),$ $\dots, \Delta_d(t;N))$ at a fixed time moment $t > 0$. In Section 5, some statistical issues having the most direct bearing on biological applications are discussed. Finally, some concluding remarks and suggestions are given in Section 6.

**2. Biological background and motivation.** The theory of branching stochastic processes has enjoyed a long history of biological applications (see, e.g., Jagers [16], Yakovlev and Yanev [22], Kimmel and Axelrod [17], Haccou, Jagers and Vatutin [10], Yanev et al. [27], Yakovlev and Yanev [25]). This theory has been proven especially useful in cell kinetics studies, many of which call for a model with multiple types of cells. A typical example is provided by the process of oligodendrocyte generation in cell culture. Oligodendrocyte type-2 astrocyte progenitor cells, henceforth referred to as O-2A progenitor cells, are known to be precursors of oligodendrocytes in the developing central nervous system. When plated *in vitro* and stimulated to divide by purified cortical astrocytes or by platelet-derived growth factor, these cells grow in clones giving rise to oligodendrocytes. An O-2A progenitor cell is partially committed to differentiation into an oligodendrocyte but it retains the ability to proliferate. Oligodendrocytes are terminally differentiated (mature) cells and they do not divide under normal conditions. At different time points over a period of several days after plating, the composition of each clone is examined microscopically to count the numbers of O-2A progenitor cells and oligodendrocytes per clone. A certain number $N$ of cell clones, each originating from a single initiator cell, are followed-up with the observation process being either longitudinal or serial sacrifice, depending on the experimental design. Then branching process modeling is needed to estimate all the important but unobservable parameters, such as the mean mitotic cycle time, probabilities of cell division and differentiation, the mean life-span of oligodedrocytes, and so forth, from the observed counts of O-2A progenitor cells and oligodendrocytes as functions of time.

The first stochastic model of oligodendrocyte development in cell culture was proposed by Yakovlev, Mayer-Proschel and Noble [23]. The model structure was defined following a set of assumptions that specified it as a special case of the Bellman–Harris branching process with two types of cells similar to that studied by Jagers [15]. Further studies [12, 13, 23, 24] proved this model to be overly simplistic, suggesting a number of refinements that have made it much more difficult to handle analytically and numerically. In parallel, estimation techniques have been developed to fit improved versions of the model to various experimental data. Because of complexity of the underlying model, these techniques have been built on simulation counterparts of either maximum likelihood or maximum pseudo-likelihood methods. All specific applications of the proposed model have invariably been limited to statistical inference from counts of both types of cells as functions of the time elapsed after plating.

The above-described approach may not be feasible in many other experimental settings. For example, it is technically impossible to count the total number of cells of a given type in the blood or bone marrow in animal experiments. A similar obstacle may arise when studying suspension cell cultures consisting of those cell types for which no specific antibodies are available, thereby limiting the utility of flow cytometry in counting cell numbers. However, this difficulty can be overcome by analyzing the proportions of different types of cells as long as they are morphologically distinguishable. While the latter expedient is routinely practiced in experimental and clinical laboratories, mathematical models of cell population kinetics are traditionally formulated in terms of cell counts. This call for making multitype branching stochastic models suitable for statistical inference from relative frequencies, rather than absolute cell counts, motivated the present work. A useful asymptotic result reported in this paper is a first step toward pertinent inferential procedures. The reported result holds for any branching process with a finite number of cell types given that the basic postulates for the evolution of each type are all met—the multitype Bienaymé–Galton–Watson, Markov, Bellman–Harris, Sevastyanov and Crump–Mode–Jagers processes representing typical examples.

**3. Basic notions and preliminary results.** Consider the process $\mathbf{Z}(t)$ introduced in Section 1 and its probability generating function (p.g.f.) given by

$$F(t;\mathbf{s}) = \mathbb{E}\{\mathbf{s}^{\mathbf{Z}(t)} \mid Z_1(0)=1\} = \mathbb{E}\{s_1^{Z_1(t)} s_2^{Z_2(t)} \cdots s_d^{Z_d(t)} \mid Z_1(0)=1\}, \tag{1}$$

where $\mathbf{s} = (s_1, s_2, \dots, s_d)$ and $|s_k| \le 1, k = 1, 2, \dots, d$.

By the branching property and the usual assumption of independent individual cell evolutions, one has

$$F_N(t;\mathbf{s}) = \mathbb{E}\{\mathbf{s}^{\mathbf{Z}(t)} \mid Z_1(0)=N\} = F^N(t;\mathbf{s}). \tag{2}$$

It is well known that all characteristics of the process can be obtained from the p.g.f. (1). In the sequel, the following notation will be used:

$$m_k(t) = \mathbb{E}\{Z_k(t) \mid Z_1(0)=1\} = \frac{\partial}{\partial s_k} F(t;\mathbf{s})|_{\mathbf{s}=\mathbf{1}}, \qquad k = 1, 2, \dots, d, \tag{3}$$

$$b_{ij}(t) = \frac{\partial^2}{\partial s_i \, \partial s_j} F(t;\mathbf{s})|_{\mathbf{s}=\mathbf{1}}, \qquad i, j = 1, 2, \dots, d, \tag{4}$$

where $\mathbf{1} = (1, 1, \dots, 1)$. Using (3) and (4), one can obtain the following moments:

$$\sigma_k^2(t) = \mathrm{Var}\{Z_k(t)\} = b_{kk}(t) + m_k(t) - m_k^2(t), \qquad k = 1, 2, \dots, d, \tag{5}$$

$$C_{ij}(t) = \mathrm{Cov}\{Z_i(t), Z_j(t)\} = b_{ij}(t) - m_i(t)m_j(t), \qquad i \ne j = 1, 2, \dots, d. \tag{6}$$

We assume that the covariance matrix $\mathbf{C}^{(d)}(t) = \|C_{ij}(t)\|$ is finite and denote its diagonal elements by $C_{ii}(t) \equiv \sigma_i^2(t), i = 1, 2, \ldots, d$. To avoid trivially unreasonable cases, we assume that $\sigma_i^2(t) > 0, i = 1, 2, \ldots, d$. Then the correlation matrix $\mathbf{R}^{(d)}(t) = \|r_{ij}(t)\|$ is well defined, where $r_{ij}(t) = C_{ij}(t)/\sigma_i(t)\sigma_j(t) = \mathrm{Cor}(Z_i(t), Z_j(t))$ and obviously $r_{ii}(t) \equiv 1, i = 1, 2, \ldots, d$.

Introduce the total number of cells at the moment $t$ as

$$U(t) = \sum_{k=1}^{d} Z_k(t), \tag{7}$$

so that the relative frequencies (fractions, proportions) of types $\Delta_k(t)$ can be defined on the nonextinction set $\{U(t) > 0\}$ as follows:

$$\Delta_k(t) = Z_k(t)/U(t), \qquad k = 1, 2, \ldots, d, \tag{8}$$

with the obvious condition

$$\sum_{k=1}^{d} \Delta_k(t) = 1. \tag{9}$$

In what follows, we will also need the following proportions:

$$p_i(t) = m_i(t)/M(t), \qquad i = 1, 2, \ldots, d, \tag{10}$$

where

$$M(t) = \mathbb{E}\{U(t)\} = \sum_{j=1}^{d} m_j(t). \tag{11}$$

To emphasize the dependence of the process $\mathbf{Z}(t)$ on the initial number of initiator cells $Z_1(0) = N$, we will use the notation $\mathbf{Z}(t; N) = (Z_1(t; N), Z_2(t; N), \ldots, Z_d(t; N))$. By the property of independence of cell evolutions, one has

$$Z_i(t; N) = \sum_{k=1}^{N} Z_i^{(k)}(t), \qquad i = 1, 2, \ldots, d, \tag{12}$$

where $\{Z_i^{(k)}(t)\}_{k=1}^{N}$ are i.i.d. copies of the branching process $Z_i(t)$, $i = 1, 2, \ldots, d$.

Our focus is on the asymptotic (as $N \to \infty$) behavior of the fractions

$$\Delta_i(t; N) = Z_i(t; N)/U(t; N), \qquad i = 1, 2, \ldots, d, \tag{13}$$

where

$$U(t; N) = \sum_{i=1}^{d} Z_i(t; N) = \sum_{k=1}^{N} U^{(k)}(t) > 0 \tag{14}$$

and $U^{(k)}(t) = \sum_{i=1}^{d} Z_i^{(k)}(t)$.

PROPOSITION 1. *Let $m_i(t) < \infty, i = 1, 2, \ldots, d$. Then for $N \to \infty$*

$$q(t;N) = \Pr\{U(t;N) = 0\} \to 0, \tag{15}$$

$$\Delta_i(t;N) \to p_i(t) \qquad a.s.\ \mathbb{E}\{\Delta_i(t;N)\} \to p_i(t), \tag{16}$$

$$\operatorname{Var}\{\Delta_i(t;N)\} \to 0. \tag{17}$$

PROOF. Note that by (1)–(2) and (12)–(14) one has $q(t;N) = \Pr\{U(t;N) = 0\} = q^N(t) \to 0$ as $N \to \infty$, where $q(t) = \Pr\{U(t) = 0\} = F(t;\mathbf{0})$ is the extinction probability, $\mathbf{0} = (0, 0, \ldots, 0)$ is a zero-vector, and $q(t) < 1$ for every fixed $t$.

From (12)–(14) and the law of large numbers (LLN), one obtains

$$\Delta_i(t;N) = \left[\frac{1}{N}\sum_{k=1}^{N} Z_i^{(k)}(t)\right] \Big/ \left[\frac{1}{N}\sum_{k=1}^{N} U^{(k)}(t)\right] \to p_i(t) \qquad \text{a.s. } N \to \infty. \tag{18}$$

Put another way, the fractions $\Delta_i(t;N)$, $i = 1, \ldots, d$, are strongly consistent estimators for $p_i(t)$ when considered as functions of the initial number of ancestors $N$. Note that for each $i$ the quantity $p_i(t)$ may be interpreted as the probability for a randomly chosen cell at time $t$ to be of the type $T_i$.

By virtue of the fact that $0 \le \Delta_i(t;N) \le 1$ and the dominated convergence theorem (DCT), it follows from (18) that $\mathbb{E}\{\Delta_i(t;N)\}$ converges to $p_i(t)$ as $N \to \infty$, implying that $\Delta_i(t;N)$ is an asymptotically unbiased estimator for $p_i(t), i = 1, 2, \ldots, d$. Similarly, by (18), one has that $\Delta_i^2(t;N) \to p_i^2(t)$ a.s. and, therefore, $\mathbb{E}\{\Delta_i^2(t;N)\} \to p_i^2(t)$ by the DCT. Hence, the result (17) follows immediately from (16). This concludes the proof of Proposition 1. □

**4. Asymptotic multivariate normality of the relative frequencies.** In the general case with $Z_1(0) = N$, the following asymptotic results hold.

THEOREM 1. *Assume $\sigma_i^2(t) < \infty, i = 1, 2, \ldots, d$, and define $A^{(d)}(t) = \|a_{ij}(t)\|$, where*

$$a_{ii}(t) = \sigma_i(t)(1 - p_i(t)), \qquad a_{ij}(t) = -\sigma_i(t)p_j(t) \quad \text{for } i \ne j;\ i, j = 1, 2, \ldots, d. \tag{19}$$

*Then for the r.v.*

$$W_i(t;N) = M(t)\sqrt{N}[\Delta_i(t;N) - p_i(t)], \qquad i = 1, 2, \ldots, d, \tag{20}$$

*the following statements are valid as $N \to \infty$:*

(i) *For every $i = 1, 2, \ldots, d$,*

$$W_i(t;N) \xrightarrow{d} Y_i(t), \tag{21}$$

*where $Y_i(t)$ is a normally distributed r.v. with $\mathbb{E}\{Y_i(t)\} = 0$ and*

$$S_i^2(t) = \operatorname{Var}\{Y_i(t)\} = \sum_{k,l=1}^{d} r_{kl}(t) a_{ki}(t) a_{li}(t). \tag{22}$$

(ii) *For every $k = 2, 3, \dots, d-1$ and every subset $(\alpha_1, \alpha_2, \dots, \alpha_k)$ with nonrecurring elements from the set $\{1, 2, \dots, d\}$,*

$$(W_{\alpha_1}(t; N), \dots, W_{\alpha_k}(t; N)) \xrightarrow{d} (Y_{\alpha_1}(t), \dots, Y_{\alpha_k}(t)), \tag{23}$$

*where $(Y_{\alpha_1}(t), \dots, Y_{\alpha_k}(t))$ have a joint normal distribution.*

(iii) *The covariance matrix of the vector $Y^{(k)}(t) = (Y_1(t), Y_2(t), \dots, Y_k(t))$ is given by*

$$D^{(k)}(t) = \| \operatorname{Cov}\{Y_i(t), Y_j(t)\} \| = [A_{d\times k}(t)]^T R^{(d)}(t) A_{d\times k}(t), \tag{24}$$

*where $A_{d\times k}(t) = \|a_{ij}(t)\|, i = 1, 2, \dots, d; j = 1, 2, \dots, k$, is a $[d \times k]$-submatrix of $A^{(d)}(t)$ and $[A_{d\times k}(t)]^T = \|a_{ji}(t)\|, j = 1, 2, \dots, k; i = 1, 2, \dots, d$, is the corresponding transposed matrix of $[k \times d]$ dimensions. The covariance matrix of any subvector $(Y_{\alpha_1}(t), \dots, Y_{\alpha_k}(t))$ can be obtained in a form similar to (24).*

PROOF. From (12)–(14), it follows that for every $i = 1, 2, \dots, d$

$$\begin{aligned} &\Delta_i(t; N) - p_i(t) \\ &\quad = \frac{\sqrt{N}}{U(t; N)} \left\{ \sigma_i(t)[1 - p_i(t)] V_i(t; N) - p_i(t) \sum_{j \neq i}^{d} \sigma_j(t) V_j(t; N) \right\}, \end{aligned} \tag{25}$$

where $V_i(t; N) = \sum_{k=1}^{N} [Z_i^{(k)}(t) - m_i(t)] / [\sigma_i(t)\sqrt{N}], i = 1, 2, \dots, d$.

Note that $\mathbb{E}\{V_i(t; N)\} = 0, \operatorname{Var}\{V_i(t; N)\} = 1$ and

$$\operatorname{Cor}\{V_i(t; N), V_j(t; N)\} = \operatorname{Cor}\{Z_i(t), Z_j(t)\} = r_{ij}(t) = C_{ij}(t)/\sigma_i(t)\sigma_j(t).$$

Then by the CLT for i.i.d. vectors (see, e.g., [2]), one has

$$(V_1(t; N), \dots, V_d(t; N)) \xrightarrow{d} (X_1(t), \dots, X_d(t)), \qquad N \to \infty, \tag{26}$$

where the r.v.s. $\mathbf{X}^{(d)}(t) = (X_1(t), \dots, X_d(t))$ have a joint normal distribution with

$$\mathbb{E}\{X_i(t)\} = 0, \qquad \operatorname{Var}\{X_i(t)\} = 1,$$
$$\operatorname{Cov}\{X_i(t), X_j(t)\} = r_{ij}(t) = C_{ij}(t)/\sigma_i(t)\sigma_j(t).$$

One can now infer from (20), (25) and (26) that the following convergence in distribution holds:

$$W_i(t;N) \xrightarrow{d} Y_i(t) = \sigma_i(t)[1-p_i(t)]X_i(t) - p_i(t)\sum_{k\neq i}^{d}\sigma_k(t)X_k(t), \qquad N\to\infty, \tag{27}$$

observing the fact that $U(t;N)/N \to M(t)$ a.s. in accordance with the LLN.

From (27) and (19), it follows that for every $i = 1,2,\ldots,d$

$$Y_i(t) = \sum_{k=1}^{d} a_{ki}(t)X_k(t) \tag{28}$$

is a linear combination of multivariate normal r.v.s. so that $Y_i(t)$ is normally distributed (see, e.g., [8], Chapter 3). Then from (28), one has

$$\mathrm{Var}\{Y_i(t)\} = \sum_{k=1}^{d}\sum_{l=1}^{d}\mathbb{E}\{X_k(t)X_l(t)\}a_{ki}(t)a_{li}(t).$$

Using (19) one arrives at (22). Formula (24) follows directly from (26) and (27). On the other hand, one can use (28) to write

$$\mathbf{Y}^{(k)}(t) = \mathbf{X}^{(d)}(t)\mathbf{A}_{d\times k}(t), \qquad k = 1,2,\ldots,d. \tag{29}$$

Formula (24) follows from (29) since the vector $\mathbf{Y}^{(k)}(t)$ is a linear transformation of the multivariate normal vector $\mathbf{X}^{(d)}(t)$ with covariance matrix $\mathbf{R}^{(d)}(t) = \|r_{ij}(t)\|$ . $\square$

Corollary 1. *By (19) and (22), it is not difficult to derive*

$$\begin{aligned} d_{ii}(t) &\equiv S_i^2(t) \\ &= (1-p_i(t))^2\sigma_i^2(t) + p_i^2(t)\sum_{k,l\neq i}C_{kl}(t) - 2p_i(t)(1-p_i(t))\sum_{k\neq i}C_{ik}(t). \end{aligned} \tag{30}$$

*To prove this, one only has to check that*

$$r_{kl}(t)a_{ki}(t)a_{li}(t) = \begin{cases} \sigma_i^2(t)(1-p_i(t))^2, & k=l=i, \\ -p_i(t)(1-p_i(t))C_{ik}(t), & k\neq i, l=i, \\ -p_i(t)(1-p_i(t))C_{il}(t), & k=i, l\neq i. \end{cases}$$

*While condition (9) implies that the fractions $\Delta_i(t;N), i = 1,2,\ldots,d,$ are linearly dependent, there exist $d-1$ joint normal distributions of lower dimensions that are asymptotically nondegenerate. For example, consider the vector*

$$\mathbf{\Delta}^{(d-1)}(t;N) = (\Delta_1(t;N),\Delta_2(t;N),\ldots,\Delta_{d-1}(t;N)).$$

*This vector has a limiting $(d-1)$-dimensional joint normal distribution with mean-vector $\mathbb{E}\{\boldsymbol{\Delta}^{(d-1)}(t;N)\} = (p_1(t), p_2(t), \ldots, p_{d-1}(t))$ and covariance matrix $\mathbf{D}^{(d-1)}(t)/NM^2(t)$, where the matrix $\mathbf{D}^{(d-1)}(t) = \|d_{ij}(t)\|$ is defined by (24) with $k = d-1$.*

COROLLARY 2. *If $i \neq j$, then*

$$d_{ij}(t) = C_{ij}(t) + p_i(t)p_j(t)\sum_{k,l} C_{kl}(t) - p_i(t)\sum_k C_{kj}(t) - p_j(t)\sum_l C_{il}(t). \tag{31}$$

*This result follows immediately from (19) and (24). Indeed, from (24), one has*

$$d_{ij}(t) = \sum_{k,l=1}^{d} r_{kl}(t)a_{ki}(t)a_{lj}(t),$$

*and using (19) one obtains*

$$r_{kl}(t)a_{ki}(t)a_{lj}(t) = \begin{cases} p_i(t)p_j(t)C_{kl}(t), & k \neq i, l \neq j, \\ -p_i(t)(1-p_j(t))C_{kj}(t), & k \neq i, l = j, \\ -(1-p_i(t))p_j(t)C_{il}(t), & k = i, l \neq j, \\ (1-p_i(t))(1-p_j(t))C_{ij}(t), & k = i, l = j. \end{cases}$$

**5. Statistical applications in relation to cell proliferation.** As was pointed out in the Introduction and in Section 2, there are experimental situations where analyzing the relative frequencies, $\Delta_i(t;N)$, of cell types rather than the total cell counts $Z_i(t;N), i = 1, 2, \ldots, d$ may be quite advantageous. Should this be the case, the property of asymptotic normality could be useful in developing the needed statistical inference of model parameters from experimental data. In particular, the following observation process is directly relevant to quantitative studies of proliferation, differentiation and death of cells. Suppose that the process under study begins with $N = \sum_{k=1}^{n} N_k$ cells of type $T_1$ and the values of $N_k$ are all large, that is, $N_0 = \min\{N_1, N_2, \ldots, N_n\} \to \infty$. The descendants of the first $N_1$ ancestors are examined only once at time $t_1$ to determine empirical counterparts of $\Delta_i(t_1;N_1), i = 1, 2, \ldots, d$, whereupon the observation process is discontinued (i.e., the cells under examination are destroyed). At the next moment $t_2 \geq t_1$, the fractions $\Delta_i(t_2;N_2), i = 1, 2, \ldots, d$, related to the descendants of the second $N_2$ ancestors are observed, and so on. This procedure results in $n$ independent observations of the form:

$$\boldsymbol{\zeta}_k = \boldsymbol{\Delta}(t_k;N_k) = (\Delta_1(t_k;N_k), \Delta_2(t_k;N_k), \ldots, \Delta_d(t_k;N_k)),$$
$$t_1 \leq t_2 \leq \cdots \leq t_n,\ k = 1, 2, \ldots, n,$$

with each vector $\boldsymbol{\zeta}_k$ being asymptotically normal in accordance with Theorem 1. Denoting the corresponding contribution to the asymptotic log-likelihood function by $L_k(\boldsymbol{\zeta}_k; t_k, N_k)$, the overall asymptotic log-likelihood is given by

$$\Lambda_n(\boldsymbol{\zeta}_1, \boldsymbol{\zeta}_2, \ldots, \boldsymbol{\zeta}_n) = \sum_{k=1}^{n} L_k(\boldsymbol{\zeta}_k; t_k, N_k). \tag{32}$$

The asymptotic log-likelihood (32) depends solely on parameters of individual multitype processes arising from a single initiator cell of type $T_1$. It is the latter parameters that are of primary interest in cell kinetics studies; they can be estimated from the data on relative frequencies by maximizing the log-likelihood (32). It should be emphasized that the only rationale for resorting to the asymptotic likelihood is that the ordinary likelihood is not readily available for partially observed branching processes of such complexity. Two more specific examples are given below.

EXAMPLE 1. Let $d = 2$ and assume that $\mathbf{Z}(t; N) = (Z_1(t; N), Z_2(t; N))$, $t = 0, 1, 2, \ldots,$ is a Bienaymé–Galton–Watson (BGW) branching process. This model represents a powerful tool in the analysis of time-lapse data generated via video-recording of individual cell evolutions. It is well known that the BGW process is entirely determined by the offspring p.g.f.

$$h_i(s_1, s_2) = \mathbb{E}\{s_1^{Z_1(1)} s_2^{Z_2(1)} \mid Z_i(0) = 1\}, \qquad i = 1, 2. \tag{33}$$

The first and second moments of the offspring distribution are derived from (32) in the usual way, that is,

$$m_{ij} = \frac{\partial}{\partial s_j} h_i(s_1, s_2)|_{s_1 = s_2 = 1}, \qquad i, j = 1, 2, \tag{34}$$

$$b^i_{jk} = \frac{\partial^2}{\partial s_j \, \partial s_k} h_i(s_1, s_2)|_{s_1 = s_2 = 1}, \qquad i, j, k = 1, 2. \tag{35}$$

Using (36) and (37), let us derive the moments defined in (3)–(6) for every generation $t = 1, 2, \ldots$. First of all, note that $\mathbf{M}(t) = \|m_{ij}(t)\| = \mathbf{M}^t$, where $m_{ij}(t) = \mathbb{E}\{Z_j(t) \mid Z_i(0) = 1\}$ and $\mathbf{M} = \|m_{ij}\|$. Let $b^i_{jk}(t) = \mathbb{E}\{Z_j(t)(Z_k(t) - \delta_{jk}) \mid Z_i(0) = 1\}$, where $\delta_{jk} = 1, j = k$ and $\delta_{jk} = 0, j \neq k$. Note also that $m_{ij}(1) = m_{ij}$ and $b^i_{jk}(1) = b^i_{jk}$.

Using the recurrence formula,

$$b^i_{jk}(t+1) = \sum_{l=1}^{2} \sum_{r=1}^{2} b^i_{lr} m_{lj}(t) m_{rk}(t) + \sum_{l=1}^{2} m_{il} b^l_{jk}(t), \qquad t = 1, 2, \ldots, \tag{36}$$

it follows from Theorem 1 and formula (32) that

$$\Lambda_n(\zeta_1,\zeta_2,\dots,\zeta_n) = -\frac{n}{2}\log 2\pi - \frac{1}{2}\sum_{k=1}^{n}\log S^2(t_k;N_k) - \frac{1}{2}\sum_{k=1}^{n}[\zeta_k - p(t_k)]^2/S^2(t_k;N_k), \tag{37}$$

where $\zeta_k = \Delta(t_k;N_k), p(t_k) = m_{11}(t_k)/M(t)$ and $M(t) = [m_{11}(t_k) + m_{12}(t_k)]$.

Taking (30) into account, one obtains

$$S^2(t_k;N_k) = \frac{1}{N_k M^2(t_k)}\{\sigma_1^2(t_k)[1-p(t_k)]^2 + \sigma_2^2(t_k)p^2(t_k) - 2C_{12}(t_k)p(t_k)[1-p(t_k)]\},$$

where, by virtue of (5) and (6), one has to set $\sigma_i^2(t_k) = b_{ii}^1(t_k) + m_{ii}(t_k) - m_{ii}^2(t_k), i = 1,2,$ and $C_{12}(t) = b_{12}^1(t) - m_{11}(t)m_{12}(t)$.

The log-likelihood $\Lambda_n$ can now be constructed [proceeding from (37)] as a function of the observations $(\zeta_1,\zeta_2,\dots,\zeta_n)$ and the moments $\{m_{ij}\}$ and $\{b_{jk}^i\}$, $i,j,k = 1,2$, to obtain asymptotic maximum-likelihood estimates of the model parameters of interest.

Example 2. Recalling the model of oligodendrocyte generation in cell culture discussed in Section 2, consider a two-type reducible Bellman–Harris process $(Z_1(t;N), Z_2(t;N))$, $t \ge 0$, with offspring p.g.f.

$$h_1(s_1,s_2) = p_0 + p_1 s_1^2 + p_2 s_2, \qquad h(1,1) = p_0 + p_1 + p_2 = 1. \tag{38}$$

In this process, the life-span of every progenitor (type $T_1$) has cumulative distribution function $G_1(t)$. At the end of its life (mitotic cycle), every progenitor cell either dies with probability $p_0$, or divides into two new $T_1$ cells with probability $p_1$, or differentiates into a new cell type $T_2$ (oligodendrocyte) with probability $p_2$. Every cell of type $T_2$ has its life-span distributed in accordance with $G_2(t)$ and, at the end of its life, it dies without giving rise to any progeny, that is, its offspring p.g.f. is $h_2(s_1,s_2) \equiv 1$.

By conditioning on the evolution of the first progenitor cell and applying the law of the total probability, one can establish that the p.g.f.s

$$F_1(t;s_1,s_2) = \mathbb{E}\{s_1^{Z_1(t)} s_2^{Z_2(t)} \mid Z_1(0) = 1\}, \qquad F_1(0;s_1,s_2) = s_1,$$
$$F_2(t;s_2) = \mathbb{E}\{s_2^{Z_2(t)} \mid Z_2(0) = 1\}, \qquad F_2(0;s_2) = s_2,$$

satisfy the following equations:

$$F_1(t;s_1,s_2) = \int_0^t h_1(F_1(t-u;s_1,s_2), F_2(t-u;s_2))\,dG_1(u) + s_1(1-G_1(t)),$$
$$F_2(t;s_2) = s_2(1-G_2(t)) + G_2(t).$$

The first and second moments of the process are readily derived from these equations. Then Theorem 1 can be applied to obtain maximum likelihood estimates of the parameters incorporated into the model.

The above line of reasoning applies to the so-called clonal analysis, that is, the analysis of cell counts in mixed clones obtained at discrete moments of time. However, in time-lapse experiments (see, e.g., [14]), it is also possible to observe the fractions $\Delta_i(t_k)$ in the embedded discrete time branching structure determined by the sizes of consecutive generations. In general, the embedded discrete time process of a $d$-type Bellman–Harris branching process is a $d$-type BGW process with the same offspring distributions. In this case, formulas (33)–(35) and (38) yield

$$\begin{aligned} &m_{11} = 2p_1, \qquad m_{12} = p_2, \qquad m_{21} = m_{22} = 0,\\ &b^1_{11} = 2p_1, \qquad b^i_{jk} = 0 \qquad \text{for all other indices } j, k,\\ &\sigma_1^2 = 4p_1(1-p_1), \qquad \sigma_2^2 = 0, \qquad C_{12} = -2p_1p_2. \end{aligned}$$

Hence, one has

$$\begin{aligned} &m_{11}(t) = (2p_1)^t, \qquad m_{12}(t) = (2p_1)^{t-1}p_2, \qquad b^1_{11}(t) = t(2p_1)^t,\\ &\sigma_1^2(t) = (2p_1)^t(t+1-2p_1), \qquad \sigma_2^2(t) \equiv 0, \qquad C_{12}(t) = -p_2(2p_1)^{2t-1} \end{aligned}$$

for $t = 1, 2, \ldots.$

Surprisingly, in this particular case

$$p(t_k) = m_{11}(t_k)/[m_{11}(t_k) + m_{12}(t_k)] \equiv 2p_1/(2p_1+p_2) \tag{39}$$

is a constant for every $t$.

Therefore, from (30) one obtains

$$\begin{aligned} S^2(t_k; N_k) &= S_1^2(t_k)/N_k M^2(t_k)\\ &= \frac{1}{N_k}\left(\frac{2p_1}{2p_1+p_2}\right)^2 \left\{t_k + 1 - 2p_1 + \left(\frac{p_2}{2p_1+p_2}\right)^2\right\}. \end{aligned} \tag{40}$$

The required log-likelihood $\Lambda_n(p_1, p_2) = \Lambda_n(\boldsymbol{\zeta}_1, \boldsymbol{\zeta}_2, \ldots, \boldsymbol{\zeta}_n \mid p_1, p_2)$ follows immediately from (39), (40) and (37).

Remark. Quine [21] investigated the moment structure of the multitype Galton–Watson process and derived useful linear recurrence relations.

**6. Discussion and concluding remarks.** Since the paper by Jagers [15], the concept of relative frequencies (cell fractions) has attracted little attention of investigators in the field of branching stochastic processes. This is unfortunate in view of the need for pertinent methods of stochastic modeling and statistical analysis of the fractions of cells rather than their counts in the area of cell biology. Published in 1971, a subsequent work by Mode [20] is evidence in favor of this opinion. Mode built his cell cycle analysis on a model of multitype positively regular age-dependent branching process. In the supercritical case, he proved that $\lim \Delta_k(t) = \delta_k$ a.s. as $t \to \infty$, providing the population does not become extinct. It should be noted that the constants $\delta_k, k = 1, 2, \ldots, d$, depend only on the offspring characteristics. In fact, $\delta_k = v_k / \sum_{j=1}^{d} v_j$, where $v_k = \eta_k(1 - G_k^*(\alpha))$, $\boldsymbol{\eta} = (\eta_1, \eta_2, \ldots, \eta_d)$ is a left eigenvector of the matrix $\mathbf{H}(\alpha) = \|G_k^*(\alpha) m_{ij}\|$ with the Malthusian parameter $\alpha$, while $G_k^*(\lambda)$ is the Laplace–Stieltjes transform of the life-span distribution $G_k(t)$ for the type $T_k, k = 1, 2, \ldots, d$. In his monograph, Mode [19] also considered the utility of fractions and reported a similar result for the BGW process.

Methods of statistical inference for branching processes with an increasing number of ancestors were developed by Yanev [26] (see also [9]), Dion and Yanev [4, 5, 6] and reviewed later by Yanev [28]. A diffusion approximation for the classical BGW process with a large number of ancestors in the near-critical case was introduced by Feller [7] and developed further by Lamperti [18] and others. Some of these results were summarized and discussed by Jagers [16]. The work of Lamperti [18] also reports some interesting limiting distributions.

However, the main focus has always been on the numbers of individuals (cells) of different types and not on their relative frequencies. In the present paper, we make another step in the same direction by considering the asymptotic behavior of the fractions $\Delta_k(t)$ as the initial number of ancestors tends to infinity but the time $t$ is fixed. The convergence results established for $\Delta_k(t)$ may have far-reaching statistical implications.

The results obtained by Jagers and Mode suggest that it would be interesting to investigate the asymptotic behavior of the fractions $\Delta_k(t; N)$ when both parameters $N$ and $t$ tend to infinity simultaneously. It is anticipated that such asymptotic properties will depend on a specific branching model and its reducibility. They are also expected to be different for supercritical, critical, and subcritical processes.

Yet another open problem has to do with correlations between times to division for sister cells. As conjectured by Harris [11], the mean number of cells is not affected by this type of correlation while the variance can only be larger than that in the independent case. Crump and Mode [3] were the first to systematically study an age-dependent branching model under which the life-spans of sister cells are correlated as well as the numbers of offspring of sister cells, but otherwise the cells live and reproduce independently. A bifurcating autoregressive branching process [10] represents another relevant example. Asymptotic properties of such processes as $N$ or/and $t$ tend to infinity have yet to be explored.

**Acknowledgments.** The authors are very thankful to the referees for their helpful comments and suggestions. This paper was prepared while N. Yanev was a visiting professor in the Department of Biostatistics and Computational Biology, University of Rochester, and he is grateful for hospitality and inspiring communications with his colleagues.

DEPARTMENT OF BIOSTATISTICS
AND COMPUTATIONAL BIOLOGY
UNIVERSITY OF ROCHESTER
601 ELMWOOD AVENUE
BOX 630
ROCHESTER, NEW YORK 14642
USA
URL: http://www.urmc.rochester.edu/smd/biostat/Andrei_Yakovlev.pdf
http://www.biology-direct.com/content/3/1/10

DEPARTMENT OF PROBABILITY AND STATISTICS
INSTITUTE OF MATHEMATICS AND INFORMATICS
BULGARIAN ACADEMY OF SCIENCES
8 G. BONCHEV
SOFIA 1113
BULGARIA
E-MAIL: yanev@math.bas.bg
Nikolay_Yanev@urmc.rochester.edu